\documentclass[notitlepage]{article}
\usepackage{graphicx}
\usepackage{amsmath}

\newtheorem{theorem}{Theorem}

\newtheorem{algorithm}[theorem]{Algorithm}

\newtheorem{corollary}[theorem]{Corollary}

\newtheorem{remark}[theorem]{Remark}

\newenvironment{proof}[1][Proof]{\textbf{#1.} }{\ \rule{0.5em}{0.5em}}
\input{tcilatex}

\begin{document}

\author{Paolo d'Alessandro \\
Dept. Of Mathematics\\
Third University of Rome\\
e-mail:dalex@math.uniroma3.it}
\title{A new conical internal evolutive LP algorithm}
\maketitle

\section{Introduction}

In this paper we develop extensions of the conical techniques exposed in the
book \cite{dacoap97}, and we develop new algorithms as well. We make an
effort to keep this paper enough self-contained, however, the book \cite
{dacoap97} (hereafter often referred to as ''the book'' for the sake of
brevity) can be useful for further details. Such book is an outgrow of two
papers published on the journal Optimization, namely\ \cite{dadualcon} and 
\cite{daprimcon92}.

To begin with the conical approach has his own nomenclature. Terms like
duality have a different meaning than those of standard approaches (see e.g. 
\cite{bertsimas1997}). Other terms like strict tangency are peculiar to the
conical approach. The term bounded refers to a LP problem (maximization of a
linear functional on a polyhedron) that is feasible and has maximum. The
specification internal for a primal algorithm has also a meaning related to
the conical context, and, as explained below, refers to the particular
approach followed here to reach optimality, as opposed to the primal
algorithm presented in the book.

In the book we stressed the particular interest of conical conditions that
are at one time purely conical ( that is, involve only cones) and purely
pointed (that is, the involved cones are pointed). Dual conical feasibility
and optimality conditions are of this kind. Although they have lead in the
book to an enumerative algorithm, the dual conical methodology is not only
useful for its peculiar properties (e.g. \ it gives a closed form expression
of the maximum,solves a whole class of problems and is particularly
efficient in handling parameter variations), but it has provided an useful
tool for applications in fields like control and game theory. The related
literature is by now very extended and relevant, and has provided many
interesting new results (see e.g \cite{dade01} and the bibliography therein).

As to the primal conical conditions given in the book, the first is conical
but not pointed. The second involves a pointed cone and an affine space, and
has been the starting point for the development an evolutive primal external
conical algorithm, enjoying exact finite convergence.

We shall take here this whole business to a more advanced stage. First
because we shall give a further primal condition, that is purely conical and
purely pointed. Secondly because, intertwining this new primal approach with
the generators machinery used in the book to develop the dual conical
algorithm, we provide a further algorithm, which although based on
generators is evolutive. Thirdly we give a closed form expression of the
maximum in the primal setting too. More than that we give explicit closed
forms of the solution both on the range space of the coefficient matrix and
in terms of unknowns of the LP problem. Finally we stress that the present
result complete the theory presented in the book from yet another point of
view. In fact the second primal optimality condition is a tangency condition
(of an affine space to the non-negative orthant) and the previous primal
algorithm reaches tangency landing the affine onto the cone. Thus we knew
that one could in principle try to reach tangency from the other side: that
is, starting with the two sets that meet each other and then taking the
affine to touch the only boundary of the cone (thereby achieving tangency).
However, developing such an algorithm (emerging submarine instead of landing
hydroplane metaphor) has defeated us for a while. It is finally presented
here.

As to notations we stick to those adopted in the book with only one
variation in the interest of simplicity. The variation is that when we shall
get to the parameterized feasibility formulation of optimality we shall
denote the augmented coefficient matrix and bound vector in the same way as
the non-augmented case, leaving to the context the specification of which is
which. Recall also that we use always the same symbol $P$ for the
non-negative orthant of any euclidean space. Again the space to which such
symbol is referred is specified by the context.

\section{Mathematical base}

We will assume thorough the hypothesis of strict tangency of the feasibility
and linear programming problems we study. Whereas this is no restriction of
generality as proved in the above book, it has a great geometrical
importance and it yields a much more neat and elegant path to the
development of algorithms. Note that actually such an assumption is made in
the primal conical LP algorithm introduced in \cite{daprimcon92}, and
reported in the book, within the theorem that states exact finite
convergence. Here however, we simplify matters avoiding to postpone the use
of such hypothesis as much as possible in an effort to achieve maximum
generality.

Let us start considering a linear feasibility problem. That is, the problem
of determining whether or not a polyhedron is non-void and in the positive
case finding at least a point in the polyhedron. Better yet would be finding
more points, i.e., having the possibility of exploring the polyhedron, even,
ideally, finding all of its points. A polyhedron is the intersection of a
finite set of semispaces. Thus a set \ $\mathcal{G}$ of the form:

\begin{equation*}
\mathcal{G}=\{x:Gx\leq \mathbf{v}\}
\end{equation*}

where $G$ is an $n\times m$ real matrix, $x$ (the unknown vector) is in $%
R^{m}$ and $\mathbf{v}$ (the bound vector)\ is in $R^{n}$. The polyhedron is
a cone if and only if \ $\mathbf{v}=0$.We shall often denote $\mathcal{R}(G)$
(the range of $G$) by $F$, for brevity. Recall also that the vector $\mathbf{%
v}-Gx$ is called the slack vector (a feasible slack vector if $\mathcal{G}%
\neq \phi $ and $x\in \mathcal{G}$).

We now make a simple but important remark, independent of the hypothesis of
strict tangency. If we decompose $\mathbf{v}$ as $\mathbf{v}=\mathbf{v}%
_{F}+\upsilon $, where $\mathbf{v}_{F}$ is the orthogonal projection of $%
\mathbf{v}$ on $F$ and $\upsilon $ is the orthogonal projection of $\mathbf{v%
}$ on $F^{^{\perp }}$, we can write for the inequality defining the
polyhedron:

\begin{equation*}
Gx\leq \mathbf{v}_{F}+\upsilon
\end{equation*}

Because vector inequalities are translation invariant, assuming $Gz=\mathbf{v%
}_{F}$, the latter is true if and only if:

\begin{equation*}
Gx-\mathbf{v}_{F}=G(x-z)\leq \upsilon
\end{equation*}

Thus in passing from one to the other condition, feasibility is invariant,
the slack vector too, and the solution is varied by a constant vector $z$.
Note that the polyhedron is a cone (when $\mathbf{v}_{F}=0$) or a translated
cone (when $\mathbf{v}_{F}\neq 0$). if and only if $\upsilon =0$. In what
follows we assume that this is not the case i.e.: $\upsilon \neq 0$.
Otherwise we always have the trivial solution $x=0$ or, respectively, $x=z$,
and our investigation would become pointless. Also notice that under this
assumption the slack vector cannot be zero because $\upsilon \neq
0\Longrightarrow \forall w\in R^{m},\upsilon -Gw\neq 0$.

\begin{remark}
Notice that if $\mathbf{v\in P}$, then $\mathcal{G}\neq \phi $, because $x=0$
is obviously a solution. Similarly if $\upsilon \mathbf{\in }P$, then $%
\mathcal{G}\neq \phi $, because $x=z$ is obviously a solution.
\end{remark}

To say that the associated problem is strictly tangent \cite{dacoap97} means
that the subspace $\mathcal{R}(G)$ is strictly tangent to the non-negative
orthant $P$, that is, $\mathcal{R}(G)$ meets $P$ in the only origin:

\begin{equation*}
\mathcal{R}(G)\cap P=\{0\}
\end{equation*}

We now recall the second primal \ conical feasibility (necessary and
sufficient) condition \cite{dacoap97}:

\begin{equation*}
\mathcal{G}\neq \phi \Longleftrightarrow (\mathbf{v}+\mathcal{R}(G)\cap
P\neq \phi
\end{equation*}

This is a conical condition, but there is an affine set appearing in it,
namely $\mathbf{v}+\mathcal{R}(G)$.

Recall also that, if feasibility prevails, to derive a feasible slack vector
in the range space of \ $G$ we can find any vector $y$ in the intersection:

\begin{equation*}
y\in (\mathbf{v}+\mathcal{R}(G)\cap P
\end{equation*}

Consequently, if we want a solution in the domain space, we can solve in $x$
the system:

\begin{equation*}
Gx=\mathbf{v}-y
\end{equation*}

Obviously if we solve this equation for all $y$ in $(\mathbf{v}+\mathcal{R}%
(G)\cap P$ we get all and nothing but the points of the polyhedron $\mathcal{%
G}$.

\section{New primal conical feasibility conditions}

We shall now start elaborating further on this second primal feasibility
condition. One of the dividends will be the introduction of a new purely
conical and pointed primal condition. As a first remark note that (as it is
immediate to prove):

\begin{equation*}
\mathbf{v}+\mathcal{R}(G)=\upsilon +\mathcal{R}(G)
\end{equation*}

and therefore we can write:

\begin{equation*}
\mathcal{G}\neq \phi \Longleftrightarrow (\upsilon +\mathcal{R}(G)\cap P\neq
\phi
\end{equation*}

Next we consider a subspace that is extended with respect to $\mathcal{R}(G)$%
, namely (with $\mathcal{L}(.)$ we denote linear extensions and with $Co(.)$
conical extensions):

\begin{equation*}
F_{e}=\mathcal{L}(\mathbf{v})+F=\mathcal{L}(\upsilon )+F=\mathcal{L(}%
\upsilon +\mathcal{R}(G))
\end{equation*}

Finally we need the cone:

\begin{equation*}
C_{e}=Co(\upsilon )+F=Co(\upsilon +F)
\end{equation*}

Note that clearly:

\begin{equation*}
F_{e}=C_{e}+(-C_{e})
\end{equation*}

We are now in a position that allows to state the following new primal
conical feasibility condition, which requires that a pointed cone do not
degenerate into the trivial (singleton of the origin) cone:

\bigskip

\begin{theorem}
The following primal conical feasibility conditions hold:

\begin{equation*}
\mathcal{G}\neq \phi \Longleftrightarrow C_{e}\text{ }\cap \text{ }P\neq
\{0\}
\end{equation*}

or equivalently:

\begin{equation*}
\mathcal{G}=\phi \Longleftrightarrow C_{e}\text{ }\cap \text{ }P=\{0\}
\end{equation*}
\end{theorem}

\begin{proof}
If $\mathcal{G}\neq \phi $, by the above version of the second primal
feasibility condition, $(\upsilon +\mathcal{R}(G)\cap P\neq \phi $ and this
intersection is made up of non-zero vectors. Let $y$ be one of those, so
that $y=\upsilon +z$ with $z$ in $\mathcal{R}(G)$. Therefore there are
non-zero points in $C_{e}$ $\cap $ $P$ and the condition is necessary.

Conversely take $y\neq 0$ in $C_{e}$ $\cap $ $P.$ It will have the
expression $y=\beta \upsilon +z$ with $z$ in $\mathcal{R}(G)$ for some real $%
\beta >0$ (if $\beta $ were zero $y$ would be zero too by virtue of strict
tangency). \ Thus the non-zero vector $w=1/\beta (\beta \upsilon +z)$
satisfies: 
\begin{equation*}
w=1/\beta (\beta \upsilon +z)=\upsilon +\frac{1}{\beta }z\in C_{e}
\end{equation*}

But evidently it is also true that $w\in (\upsilon +\mathcal{R}(G))\cap P$
and this completes the proof.
\end{proof}

\bigskip

The procedure explained in the proof to pass from the vector $y$\ in the
cone $C_{e}$ $\cap $ $P$ to the vector $w$ moving along the ray generated by 
$y$, will be called calibration. Notice that such procedure is viable
numerically in a very simple way. In fact, if we denote by $P_{F^{\bot }}$
the orthogonal projection onto $F^{\bot },$ it is obvious that:

\begin{equation*}
\beta =(P_{F^{\bot }}y)_{i}/\upsilon _{i};\forall i\text{ s.t. }\upsilon
_{i}\neq 0
\end{equation*}

Another important fact that follows immediately from the above proof is that
if we let, in the feasibility case $P_{c}=(v+\mathcal{R}(G)\cap P$ (more on
this set later) we can state the following:

\begin{corollary}
The following relation holds whenever feasibility prevails

\begin{equation*}
Co(P_{c})=C_{e}\text{ }\cap \text{ }P
\end{equation*}
\end{corollary}

Next notice that if $\mathcal{G}=\phi $ so that $C_{e}$ $\cap $ $P=\{0\}$ it
may either be the case that:

\begin{equation*}
-C_{e}\text{ }\cap \text{ }P=\{0\}
\end{equation*}

in which case it is clear that $F_{e}$ is strictly tangent to $P$; or that

\begin{equation*}
-C_{e}\text{ }\cap \text{ }P\neq \{0\}
\end{equation*}

In this respect, the next natural question is to find out whether it can be
the case that $-C_{e}$ $\cap $ $P\neq \{0\}$ and $C_{e}$ $\cap $ $P\neq
\{0\} $ at the same time. We shall show in the proof of the next result that
this cannot be the case, because the two cones $C_{e}$ and $-$ $C_{e}$ live
in opposite semispaces and the non-negative orthant, with the exception of
the origin, is entirely contained in the interior of one of them. And this
means that we can establish the following further feasibility condition,
which is more handy in that it substitutes the subspace $F_{e}$ to the cone $%
C_{e}$.

\begin{theorem}
If \ $F_{e}$ $\cap $ $P\neq \{0\}$ then the sign of $\beta $ is constant $%
\forall y\neq 0$ in \ $F_{e}\cap P$. Moreover:

\begin{equation*}
\mathcal{G}\neq \phi \Longleftrightarrow \ F_{e}\text{ }\cap \text{ }P\neq
\{0\}\text{ and }\beta >0
\end{equation*}
\end{theorem}

\begin{proof}
Suppose $F_{e}$ $\cap $ $P\neq \{0\}$ and consider an $y\neq 0$ in $F_{e}$ $%
\cap $ $P$. \ Because $F$ is strictly tangent to $P$, we know from theorem
6.2.1 in \cite{dacoap97} that $F$ is contained in an hyperplane $H$, which
is strictly tangent to $P$ and there is a vector $n$, normal to this
hyperplane, which is internal to $P$. Then for some $z\in $:$F,\gamma \neq 0$
:

\begin{equation*}
0<(n,y)=(n,(\gamma \upsilon +z))=\gamma (n,\upsilon )
\end{equation*}

Assume to fix the ideas that $\gamma >0$ so that $(n,\upsilon )>0$, $y\in
C_{e}$ and therefore $y\in $.$C_{e}$ $\cap $ $P$. In this case by the
previous condition feasibility prevails. Clearly by the same argument and
because $n$ is internal to $P$:

\begin{eqnarray*}
(n,C_{e}) &\subset &[0,+\infty ) \\
(n,P) &\subset &[0,+\infty )
\end{eqnarray*}

Thus both the cones $C_{e}$ and $P$ are contained in the semispace $%
\{x:(n,x)\geq 0\},$ delimited by the hyperplane $H$, and, in addition, $P$
is in the interior of the semispace, with the only exception of the origin,
because $n$ is in the interior of $P$.

Moreover::

\begin{equation*}
(n,-C_{e})\subset (-\infty ,0]
\end{equation*}

and this latter implies that $-C_{e}$ is in the opposite semispace (i.e. the
semispace $\{x:(n,x)\leq 0\}$) with respect to $P,$ so that $-C_{e}$ $\cap $ 
$P=\{0\}$. Therefore the sign of $\beta $ in $F_{e}$ $\cap $ $%
P=(C_{e}+(-C_{e}))\cap $ $P$ is constant.

Assuming instead that $\gamma <0$ \ and reasoning along the same lines $%
(n,\upsilon )<0$, $y\in -C_{e}$ and therefore $y\in $.$-C_{e}$ $\cap $ $P$.
Moreover in this case:

\begin{equation*}
(n,C_{e})\subset (-\infty ,0]
\end{equation*}

and therefore $C_{e}$ $\cap $ $P=\{0\}$ and feasibility does not prevail. In
this case $-C_{e}$ $\cap $ $P\neq \{0\}$ and the sign of $\beta $ is
constant in $-C_{e}$ $\cap $ $P$. Thus the proof is finished.
\end{proof}

\begin{corollary}
There are three mutually exclusive and exhaustive cases possible for the
feasibility problem:

a) $F_{e}$ is strictly tangent to $P$. That is $F_{e}$ $\cap $ $P=\{0\}$: In
this case the problem is unfeasible.

b) $F_{e}$ $\cap $ $P\neq \{0\}$ and $\beta <0$. In this case too the
problem is unfeasible.

c) $F_{e}$ $\cap $ $P\neq \{0\}$ and $\beta >0$. In this case the problem is
feasible
\end{corollary}

As in all algorithms exposed in \cite{dacoap97}, once it is determined that
the problem is feasible and an $y\in C_{e}$ $\cap $ $P$ is found, a solution
in the domain space can be found solving the equation $Gx=v-w$, where $w$ is
obtained from $y$ by the calibration procedure.

\section{Feasibility algorithm}

We next introduce a feasibility algorithm, based on the results obtained so
far. To this purpose we have to devise a method to find (in the feasible
case) one or more feasible solutions in the range space $y\in C_{e}$ $\cap $ 
$P$ . Or, which is more easily viable, find \ an $y\in F_{e}$ $\cap P$ with $%
\beta >0$ and then calibrate it. The method should fail if and only if
unfeasibility prevails so to completely solve the problem. One possible and
useful way to go is to exploit the generators machinery for pointed
polyhedral cones (and we got purely conical and pointed conditions primarily
to that effect) in this primal conical setting. The advantages of this
approach will be more and more evident in the sequel. Because $y$ must be in
a pointed polyhedral cone given by the intersection of a subspace with the
non-negative orthant, we can find one or more of the extreme rays of this
cone or, equivalently, of its generators. A good part of the book \cite
{dacoap97} is devoted to the development of a sophisticated machinery to
solve this problem, with basic algorithms and two levels of more advanced
algorithms as well as an implementation and numerical results. The
fundamental results on which the generators techniques are based are given
in Chapter 9.

To carry out this program, we only need to find an operator whose null space
is $F_{e}$. This is not difficult. In fact:

\begin{equation*}
F_{e}=\mathcal{L}(\upsilon )+F
\end{equation*}

Where the sum is a direct orthogonal sum. We introduce the notations $V=$ $%
\mathcal{L}(\upsilon )$ and $P_{V}$ and $P_{F}$ to denote, respectively, the
orthogonal projections onto $\ V$ and $F$. Then the orthogonal projection
onto $\mathcal{L}(\upsilon )+F$ is given by $P_{V}$ +$P_{F}$ so that $%
I-P_{V} $ $-P_{F}$ is the orthogonal projection onto $(V+F)^{\perp }$
Therefore:

\begin{equation*}
F_{e}=\mathcal{L}(\upsilon )+F=\mathcal{N}(I-P_{V}-P_{F})
\end{equation*}

Thus we can state the following feasibility algorithm. \ In the feasible
case one can compute either a single or more solutions in the range space.
This is not a detail it is a major point on which we shall expand in the
next section.

\begin{algorithm}
: New primal conical feasibility algorithm

1- Preliminary step: verify whether\textbf{\ }$\mathbf{v}$ or $\upsilon $
are in $P$. If either of those is the case trivial solution are immediately
available as explained in Section 2. In such a case STOP.

2- Compute $P_{V}$ and $P_{F}$ \ and hence $I-P_{V}-P_{F}$.

3- Use the techniques of \cite{dacoap97} to compute a first generator of $%
\mathcal{N}(I-P_{V}-P_{F})\cap P$. If none is found the problem is
unfeasible. Hence STOP. Otherwise go to step 4.

4- Calibrate the generator. (Recall that $\beta =0$ cannot happen) If $\beta
<0$ the problem is unfeasible. In such case STOP. If $\beta >0$ the problem
is feasible. Proceed to the next steps.

5- Facultative step. Can be executed or skipped. Go on to calculate more or
all the generators and calibrate each of them.

6- For each calibrated generator $g_{c}$ a solution can be obtained solving $%
Gx=v-g_{c}$.
\end{algorithm}

\begin{remark}
In certain applications, it may be convenient (of course in the feasible
case) to find a point in relative interior of the pointed polyhedral cone $%
F_{e}\cap P=\mathcal{N}(I-P_{V}-P_{F})\cap P$. To this effect one can
compute all the calibrated generators and then the sum of them (Theorem
6.1.2 in \cite{dacoap97}) is in the relative interior of the cone.
\end{remark}

The significance of the set of all the calibrated generator goes way beyond
this simple remark and is of crucial importance in the present context. This
issue is fully explored in the next section.

\section{The contact polytope}

Suppose that the problem is feasible. Then, as recalled above, the
polyhedron $P_{c}=(\mathbf{v}+\mathcal{R}(G)\cap P=(\upsilon +\mathcal{R}%
(G)\cap P$ is non-void. In the book, leaning on a result by Ben Israel it is
shown that this set is actually a polytope, which we called the contact
polytope. The properties of the contact polytope play an important role in
the book, but it was not determined explicitly. Here we complete the theory
with a detailed description. In fact the set $\{g_{ci}:i=1,..,k\}$ of
calibrated generators of the pointed polyhedral cone $C_{e}$ $\cap $ $P$ is
exactly the set of extreme points of $P_{c}$. This is stated in the next
theorem:

\begin{theorem}
The set of calibrated generators of the cone $C_{e}$ $\cap $ $P$ is equal to
the set of extreme points of the contact polytope $P_{c}=(\mathbf{v}+%
\mathcal{R}(G)\cap P$. Thus:

\begin{equation*}
P_{c}=\mathcal{C(\{}g_{ci}\})
\end{equation*}
\end{theorem}

\begin{proof}
First we prove that any ray of the cone $C_{e}$ $\cap $ $P$ can intersect $%
P_{c}$ only in a single point. For suppose that two points of a ray, say $%
z_{1}$ and $\ z_{2}$, be in the intersection, so that we may assume, without
restriction of generality, that $z_{2}=\alpha z_{1}$ with $\alpha >1.$ Then $%
z_{2}-z_{1}\neq 0$ would be both in $\mathcal{R}(G)$ and in $P$. But this
contradicts the hypothesis that $\mathcal{R}(G)$ is strictly tangent to $P$.
Therefore $z_{2}=z_{1}$ as we wanted to prove.

We know from the previous section that all the points of $%
\{g_{ci}:i=1,..,k\} $ are in $P_{c}$ and $P_{c}$ is contained in $C_{e}$ $%
\cap $ $P$. Let now $\{z_{j}:j=1,..,p\}$ be the set of extreme points of $%
P_{c}$.

We start showing that each \ calibrated generator is an extreme point of $%
P_{c}$. In fact suppose that is not so for a certain $\ g_{ci}$. Then we can
write (in what follow we drop indexes in sums to keep the notation simple):

\begin{equation*}
\ g_{ci}=\sum \alpha _{r}z_{r}
\end{equation*}
with at least two non-zero coefficients. Isolating one term, we can write
without restriction of generality (possibly the sum reduces to a single
term):

\begin{equation*}
\ g_{ci}=\alpha _{1}z_{1}+\sum \alpha _{\eta }z_{\eta }=\alpha
_{1}z_{1}+(\sum \alpha _{\eta })\sum (\alpha _{\eta }/\sum \alpha _{\eta
})z_{\eta }=\alpha _{1}z_{1}+\beta w
\end{equation*}
where $\alpha _{1}+\beta =1$ and $w\in P_{c}$. Because $z_{1}$ is an extreme
point, $z_{1}\neq w$, and by what we proved at the beginning, these two
vectors cannot be proportional. It follows that this equality contradicts
that $g_{ci}$ \ lies on an extreme ray and therefore we have reached the
conclusion that all the points $\ g_{ci}$ are actually extreme points of $%
P_{c}$.

Finally assume that there is some extreme points that are not in the set of
calibrated generators. Let one of these be $z_{j}$. Then because $z_{j}$.$%
\in C_{e}$ $\cap $ $P$:

\begin{equation*}
z_{j}=\sum \alpha _{i}g_{ci}
\end{equation*}

where $\alpha _{i}\geq 0$ \ $\forall $ $i$ and $\sum \alpha _{i}>0$ \ If
just one $\alpha _{i}$ is non-zero a contradiction is immediate, because if $%
\alpha _{i}=1$, then $z_{j}$ is a calibrated generator and if $\alpha
_{i}\neq 1$ there are two proportional vectors in $P_{c}$. Excluding this
case, either $\sum \alpha _{i}=1$, and then we get the contradiction that an
extreme point is a non-zero convex combination of a set of other extreme
points or $\sum \alpha _{i}\neq 1$. If that is so, dividing both sides of
the above expression of $z_{j}$ by $\sum \alpha _{i}$ we get again the
contradiction that there are two points in a ray that belong to the polytope 
$P_{c}$. The proof is thereby completed.
\end{proof}

The determination of the contact polytope and, in particular, the
identification of the extreme rays of $C_{e}$ $\cap $ $P$ with the rays
generated by the extreme points of $P_{c}$ \ is a very important addition to
our theory. We can score on this result immediately.

We start with a remark on solutions. In fact $P_{c}$ can be also viewed as
the set of all feasible slack vectors, in the sense that it contains all and
nothing but the slack vectors $y$ for which all solution of $Gx+y=\mathbf{v} 
$ is a feasible solution. Consequently the set of vectors $x$ obtained in
this way is the set of all feasible solutions. We can summarize this in the
following::

\begin{corollary}
if we let the above feasibility algorithm run to compute all the calibrated
generators, then all and nothing but the solution of the problem are given
by $\{x:Gx=\mathbf{v}-y,$ $y\in \mathcal{C(\{}g_{ci}\})\}$.
\end{corollary}

\begin{remark}
It should be stressed at this point that the present theory yields in a
primal conical setting an explicit expression for both the polytope of
feasible solution in the range space (slack vectors) and the polyhedron of
the feasible solutions of the problem
\end{remark}

Another important consequence of this result is that we are now in the
position of introducing an internal primal conical LP algorithms.

Before getting into this (in the next section), let us recall briefly the
well-known parameterized feasibility formulation on which we based all LP
methods (see the book for more details).

Consider the problem of maximizing a linear functional on a polyhedron
(Problem LP)

\begin{eqnarray*}
&&\max \text{ }fx\text{ } \\
subject\text{ }to &:&\text{ }Gx\leq \mathbf{v}
\end{eqnarray*}

We can rewrite this as:

\begin{eqnarray*}
&&\max \text{ }h \\
subject\text{ }to &:&\widehat{G}x\leq \widehat{\mathbf{v}}(h)
\end{eqnarray*}

where the augmented coefficient matrix $\widehat{G}$ is obtained adding to $%
G $ a row with the entries of $-f$ and the augmented bound vector $\widehat{%
\mathbf{v}}(h)$ \ is obtained adding to $\mathbf{v}$ \ a last entry equal to 
$-h$.

In what follows we shall soon use for this problem the same notations of the
feasibility problem in order to simplify our presentation. Thus the reader
is advised that it is the context to determine whether we refer to plain or
augmented coefficient matrix and to plain or augmented bound vector. The
same rule we apply to all the other mathematical entities related to the
problem. Thus, for example we still denote by $P_{c}$ the contact polytope
corresponding to the polyhedron defined by $\widehat{G}x\leq \widehat{%
\mathbf{v}}(h)$. Also the non-negative orthant is still denoted by $P$ in
the augmented euclidean space. However, whenever beneficial to clarity, we
explicitly denote dependence of the appropriate items on $h$ (writing, e.g., 
$P_{c}(h))$.\ 

In this parameterized feasibility setting, to solve the problem, we have to
find an $h_{o}$ such that if $h>h_{o}$ then the polyhedron $\{x:\widehat{G}%
x\leq \widehat{\mathbf{v}}(h)\}$ is void and, if $h\leq h_{o}$, \ the same
polyhedron is non-void. Such an $h_{o}$, if it exists at all, is the optimum
value of the functional. In terms of the second primal conical feasibility
condition we have to find $h_{o}$ that verifies the following tangency
condition:

\begin{equation*}
(\widehat{\mathbf{v}}(h)+\mathcal{R}(\widehat{G}))\cap P=\phi \text{ if }%
h>h_{o}
\end{equation*}

\begin{equation*}
(\widehat{\mathbf{v}}(h)+\mathcal{R}(\widehat{G}))\cap P=P_{c}\neq \phi 
\text{ if }h\leq h_{o}
\end{equation*}

The primal conical LP algorithm in the book (under conditions that insure
feasibility and boundedness - see next section) started from a large value
of $h$ , so to insure that $P_{c}$ be void and diminished such an $h$ until
the tangency of the affine space to the non-negative orthant is reached, so
that such value of $h$ is just the optimum value $h_{o}$of the functional .
This is the external approach.

Here we take the dual (how many meaning of this word!) view: we start from
low values of $h,$ so that the contact polytope is non-void, and increment $%
\ h$, until the contact polytope is squeezed into the maximal face of the
non-negative orthant, that corresponds to a zero last component. At that
point tangency prevails and, therefore, the corresponding value of $h$ is
equal to the optimum $h_{o}$. Actually we give an algorithm \ of this sort
in two versions. The second one will realize a further particularly
important advance within our methodology. To each of them we devote the next
two sections.

\section{Primal conical internal algorithm: first version.}

As in the book, and without restriction of generality, we assume strict
tangency Feasibility can be ascertained \ as illustrated in the previous
section. We recall from the book that, in view of Theorem 6.3.1, under
strict tangency, feasibility implies boundedness. Thus we can now work on
optimality with all three hypothesis (strict tangency feasibility and
boundedness) in force. Keep in mind that we refer to the augmented problem
now, although we do not change notations.

In view of the results of the foregoing section, in the feasible case, to
know the calibrated generators of $C_{e}$ $\cap $ $P$ is the same as knowing
the extreme points of the contact polytope $P_{c}$. As we shall see
momentarily the knowledge of certain extreme points of $P_{c}.$ is
equivalent to the knowledge of the maximum value of the functional . It
cannot be overemphasized the importance of the fact that, by the following
theorem, we obtain a closed form expression for the maximum in a primal
setting, just as a closed expression for the maximum was already given in
the dual setting (Theorem 5.3.1 in the book).

\begin{theorem}
\label{maxformula}Suppose to choose an $h$ \ such that the augmented problem
is feasible($h\leq h_{o}$) and let $\{g_{ci}:i=1,..,k\}$ be the set of
extreme points of $P_{c}(h)$ (or, what is the same, calibrated generators of
.$C_{e}(h)$ $\cap $ $P$) then, letting $h_{m}=\max \{g_{ci_{n+1}}:i=1,..,k\}$%
:

\begin{equation*}
h_{o}=h+h_{m}
\end{equation*}
\end{theorem}

\begin{proof}
Passing from an $h$ to $h+\Delta h$ corresponds to adding $-\Delta h$ to the
last component of all points of $\mathbf{v}+\mathcal{R}(G)$. Thus all points
\ in the contact polytope with the last component equal to $h_{m}$ are still
in the new contact polytope $P_{c}(h+h_{m})$ and fall in the maximal face of 
$P,$ $M_{n+1}=\{y:y\in P$ and \ $y_{n+1}=0\}$. \ Suppose now that in $%
P_{c}(h+h_{m})$ there is a point $z$ with $z_{n+1}=\rho >0.$ Then passing
from $h+h_{m}$ to $h$ (that is, incrementing $h$ by $-h_{m}$)\ this point $z$
translates to a point with all the same components but the last, which is
equal to $\rho +h_{m}>h_{m}$ and, in addition, it clearly is in $P_{c}(h)$.
However the last component $\gamma $ of all the points in $P_{c}(h)$ \
satisfies $h_{m}=\max \{g_{ci_{n+1}}:i=1,..,k\}$ $\geq \gamma \geq $ $\min
\{g_{ci_{n+1}}:i=1,..,k\}$ by fact that a polytope is the convex extension
of the set of its extreme points. Therefore we have found a contradiction
and it follows that $P_{c}(h+h_{m})\subset M_{n+1}$, or, in other words, $%
\mathbf{v}+\mathcal{R}(G)$ is tangent to $P$. \ \ \ By the second primal
conical optimality condition, $h_{o}=h+h_{m}$ is the optimum value of the
functional and solves the problem in the range space of $G$.
\end{proof}

In the proof of the theorem the procedure to find solutions is already built
in. We record such procedure in the following Corollary (by $e_{i}$ we
denote the vector that has all zero components but the $ith$, which is equal
to $1$):

\begin{corollary}
(i) Let \ $g_{c}$ be an extreme point of $P_{c}(h)$ such that $%
g_{c_{n+1}}=h_{m}$. Then an optimum slack vector is given by $%
y_{o}=g_{c}-h_{m}e_{n+1}$, $y_{o}$ is an extreme point of $P_{c}(h_{o})$ and
all solutions of the equation $Gx=\mathbf{v}-y_{o}$ (that necessarily exist)
are optimum solutions.

Let $\{g_{ci}\}$ be the set of extreme points of \ $P_{c}(h_{o})$ . Then the
set of all optimum solutions is obtained solving the equation $Gx=\mathbf{v}%
-y$ where $y\in \mathcal{C}(\{g_{ci}\})$.
\end{corollary}

\begin{remark}
With this Corollary we have completed the picture of the conical approach
giving also the explicit closed form for the sets of optimal solutions both
in the range space and in the domain space.
\end{remark}

The proof of the Corollary requires but trivial new verifications ans can be
safely omitted.

We can now structure a PL algorithm.. First set an $h$ small enough (it can
be arbitrarily small) to ensure that $h<h_{o}$ (a remark on this is given
right after the statement).. Then apply the following

\begin{algorithm}
(Primal conical internal algorithm)

Step1 Find all the calibrated generators with positive last component.

Step2. Set $h_{o}=h+h_{m},$ where $h_{m}$ is the maximum of last components
of the found calibrated generators

Step3 Consider any calibrated generator $g_{ci}$ such that $%
g_{ci_{n+1}}=h_{m}$. Then an optimum slack vector is given by $%
y_{o}=g_{ci}-h_{m}e_{n+1}$ and an optimum solution $x_{o}$\ is given by any
solution of the equation $Gx=\mathbf{v}-y_{o}$.
\end{algorithm}

The proof of the algorithm is given in the theory so far developed. The
above Corollary also illustrate how to find the set of all the solutions, if
needed.

\begin{remark}
Note that the internal and external algorithm (Described in Ch 11 of the
book) complete each other. In fact if no calibrated generator with positive
last component is found we are either at the optimum i.e., $h=h_{o}$, or $%
h>h_{o}$. In any case we can revert to the external algorithm \cite{dacoap97}
and find the solution.
\end{remark}

An interesting aspect of this algorithm is that it is not purely enumerative
because we do not look for all the generators. We may easily reformulate it
requiring that in STEP 1 all the generators be found. In this case the
remark is changed accordingly. If we find generators, but none has a
positive last component then $h=h_{o}.$ If no generator is found $h>h_{o}$
and we can revert to the external algorithm. \ 

We can pursue that feature further and introduce an evolutive version of the
algorithm. Although the evolutive character could be exhibited in abstract
terms, it becomes more evident if we take to the fore .the techniques
illustrated in the book to find the generators of a polyhedral pointed cone,
which is the intersection of a subspace and the non-negative orthant. We
conjugate in this way the generators technique used to deploy the dual
conical methods (which were essentially enumerative) with the present new
primal approach and fulfill the quest for evolutiveness mentioned in the
book within the dual conical framework.

\section{Computation of generators and the evolutive version of the algorithm%
}

We assume the same hypotheses and, in particular, feasibility and $h<h_{o}$
\ Our first purpose is to show that it is possible to apply to the present
problem the machinery developed in the book for the computation of
generators. The peculiarities of our method will then allow us to derive an
evolutive algorithm.

Recall that the orthogonal projection $P_{F_{e}}$ of the space onto $F_{e}=%
\mathcal{L}(\upsilon )+F=V+F,$ is given by:

\begin{equation*}
P_{F_{e}}=P_{V}+P_{F}
\end{equation*}

and the orthogonal projection of the space onto $F_{e}^{\bot }$ is given by:

\begin{equation*}
P_{F_{e}^{\bot }}=I-(P_{V}+P_{F})
\end{equation*}

It follows that we can express $F_{e}$ as:

\begin{equation*}
F_{e}=\mathcal{N}(P_{F_{e}^{\bot }})=\mathcal{N}(I-(P_{V}+P_{F}))
\end{equation*}

At this point we can apply all the machinery developed in the book to find
the generators of

\begin{equation*}
F_{e}\cap P=\mathcal{N}(I-(P_{V}+P_{F}))\cap P
\end{equation*}

Notice that in this formula only $P_{V}$ depends on $h$.

However, with respect to the case of the dual conical method, there are
numerous simplifications. First because we want those generators that have a
non-zero last component. In this respect we can state the following:

\begin{theorem}
Under the present hypotheses, \ deleting the last column of \ the matrix $%
I-(P_{\mathcal{L}(\upsilon )}+P_{F})$, we obtain a matrix with the same rank
as the original matrix.
\end{theorem}

\begin{proof}
In view of theorem 15.1.1 of the book, if it were not so, we would not get
any generator with non-zero last component, and hence in view of Theorem \ref
{maxformula} above a contradiction would arise.
\end{proof}

\begin{corollary}
In the procedure for search of generators of $\mathcal{N}(I-(P_{V}+P_{F}))%
\cap P$ with non-zero last component given in the book we can constantly use
the last column as test column.
\end{corollary}

We are now ready to introduce an evolutive conical algorithm

\begin{algorithm}
(Primal conical internal evolutive algorithm)

Put $h(0)=h$. Put $T(0)=I-(P_{V(h)}+P_{F})$. Repeat the following step for
i=1,2,..:

STEP i:\ If i=1 perform the procedure of the book modified fixing the last
column as test column to find generators of .$\mathcal{N}(T(0))\cap P.$ If i%
\TEXTsymbol{>}1 resume the search from the sequence of basic column
subsequent to the last of step i-1. If a generator is found, do not verify
it was already found, and proceed to calibrate it. Let $g(i)$ be the
calibrated generator. Set $h(i)=h(i-1)+$ $g(i)|_{n+1}.$ Compute $T(i)$
setting $h=h(i).$

Until the procedure introduced in the book terminates.

If the loop is exited at $i=j,$ set $h_{o}=h(j).$ Optimal solutions can be
obtained solving for $x$ the equation $G(x)=v(h_{o})-g(j)$.
\end{algorithm}

\begin{theorem}
The evolutive algorithm enjoys exact finite convergence. That is \ it
converges in a finite number of step and if the loop is exited at $i=j,$ $%
h_{o}=h(j).$
\end{theorem}

\begin{proof}
The proof is essentially contained in the proof of the enumerative version.
What we do is to pass from a calibrated generator to the next with
increasing last component. Because the number of calibrated generators is
finite it is granted the algorithm converges in a finite number of steps.
The only thing that remains to be proved is that the technique to find the
sequence of calibrated generators is correct. But this too is rather
obvious. In fact each time we increase $h$ the current calibrated generator
(as well as any calibrated generator with lower last component) is
eliminated. Thus any new calibrated generator has a larger last component
with respect to the former ones. and, consequently, we never have to verify
that we find already known calibrated generators. Moreover we can resume the
search from where it left at each step, because if we started from scratch
and found a calibrated generator, the same calibrated generator would have
appeared before, by the argument used in the proof of the first version of
the algorithm, and would have appeared with a lower last component. And we
know from the theory developed hitherto that this is a contradiction. The
rest of the algorithm (computation of solutions) should by now obvious.
\end{proof}

\bigskip

\begin{remark}
Notice that, as made clear by the above proof, the algorithm will usually
get rid of some of the generators that are computed in the non-evolutive
version. In other words evolutiveness is not just adjourning the value of
the maximum, but, in general, avoiding the necessity of visiting the whole
set of extreme points of the contact polytope.
\end{remark}

\section{Conclusion}

As for the previous conical algorithms we deferred submission until we had
evidence that the algorithm performs correctly numerically. During the
development of the implementation we used the same example of \cite
{dadualcon}. Of course, as in the previous cases, timeliness was
priviledged, so to arrive to a first straightforward implementation, that in
the present case was written in Pascal, within the Delphi environment..

The previous experience showed that optimization of the code is a lengthy
and painful endevour, that required the derivation of further results that
are accounted for in the book. However, we got the divident of \ entire
orders of magnitude improvings in computing time. The same process is in its
inception for the present algorithm. An uprise of its numerical efficiency
will be given as soon as we will feel that the level optimization of the
code is satisfactory, and further improvements will have marginal effects
only.

\end{document}